\documentclass[11pt,a4paper]{amsart}
\usepackage{amsfonts}
\usepackage{amssymb}
\usepackage{amsthm}
\usepackage[centertags]{amsmath}
\usepackage{newlfont}
\usepackage{graphicx}

\newtheorem{theorem}{Theorem}

\newtheorem{definition}[theorem]{Definition}

\def\r{\mathbb R}
\def\l{\mathbb L}

\begin{document}
\title{Translation and homothetical surfaces in Euclidean space with constant curvature}
\author{Rafael L\'opez}
 \address{Departamento de Geometr\'{\i}a y Topolog\'{\i}a\\
Universidad de Granada\\
18071 Granada, Spain\\}
\email{ rcamino@ugr.es}
\thanks{ The corresponding author. Partially supported by MEC-FEDER
 grant no. MTM2011-22547 and Junta de Andaluc\'{\i}a grant no. P09-FQM-5088.}
\author{Marilena Moruz}
\address{Department of Mathematics\\
A. I. Cuza University, Iasi, Romania }
\email{marilena.moruz@gmail.com}
\thanks{Partially supported by CNCS-UEFISCDI grant PN-II-RU-TE-2011-3-0017 (Romania)}

\begin{abstract} We study surfaces in Euclidean space constructed by the sum of two curves or that are graphs of the product of two functions. We consider the problem to determine all these surfaces with constant Gauss curvature. We extend the results to non degenerate surfaces in Lorentz-Minkowski space.
\end{abstract}
\subjclass[2000]{ 53A10, 53C45}
\keywords{translation surface, homothetical surface, mean curvature, Gauss curvature}

 \maketitle

\section{Introduction}

In this paper we study two types of surfaces in Euclidean space $\r^3$. The first kind of surfaces are translation surfaces which were initially introduced by S. Lie. A \emph{translation surface} $S$ is a surface that can be expressed as the sum of two curves $\alpha:I\subset\r\rightarrow\r^3$, $\beta:J\subset\r\rightarrow\r^3$. In  parametric form, the surface $S$ writes as $X(s,t)=\alpha(s)+\beta(t)$, $s\in I, t\in J$. See \cite[p. 138]{da}. A translation surface $S$ has the property that the translations of a parametric curve $s=ct$ by $\beta(t)$ remains in $S$ (similarly for the parametric curves $t=ct$).
It is an open problem to classify all   translations surfaces with constant mean curvature (CMC) or constant Gauss curvature (CGC). A first example of CMC translation surface is the Scherk surface
$$z(x,y)=\frac{1}{a}\log\left(\left|\frac{\cos(ay)}{\cos(ax)}\right|\right),\ a>0.$$
This surface is minimal ($H=0$) and it belongs to a more general family of Scherk surfaces (\cite[pp. 67-73]{ni}). In this case, $\alpha$ and $\beta$ lie in two orthogonal planes and after a change of coordinates, the surface is locally described as the graph of $z=f(x)+g(y)$. Other examples of CMC or CGC translation surfaces given as a graph $z=f(x)+g(y)$  are: planes ($H=K=0$), circular cylinders ($H=ct\not=0$, $K=0$) and cylindrical surfaces ($K=0$). The progress on this problem has been as follows.
\begin{enumerate}
\item If $\alpha$ and $\beta$ lie in orthogonal planes, the only minimal  translation surface is the Scherk surface and the plane \cite{sc}.
 \item If $\alpha$ and $\beta$ lie in orthogonal planes, the only CMC translation surface is the plane, the Scherk surface and the circular cylinder \cite{li}.
 \item If $\alpha$ and $\beta$ lie in orthogonal planes, the only CGC translation surfaces have $K=0$ and are cylindrical surfaces \cite{li}.
 \item If both curves $\alpha$ and $\beta$ lie in planes, the only minimal translation surface is a plane or it  belongs to the family of Scherk surfaces \cite {dwvw}.
 \item If one of the curves $\alpha$ or $\beta$ is planar and the other one is not, there are no minimal translation surfaces \cite{dwvw}.
 \end{enumerate}

  Our first  result concerns the case that the Gauss curvature $K$ is constant. We prove that, without any assumption on the curves $\alpha$ and $\beta$,  the  only known flat ($K=0$) translation surfaces are cylindrical surfaces. By a cylindrical surface we mean a ruled surface whose rulings are all parallel to a fix direction.

 \begin{theorem}\label{tgt}
  \begin{enumerate}
 \item The only translation surfaces with constant Gauss curvature $K=0$ are cylindrical surfaces.
 \item There are not translation surfaces with constant Gauss curvature $K\not=0$ if one of the generating curves is planar.
 \end{enumerate}
 \end{theorem}
We observe that when $K=0$ we give a complete classification CGC of translation surfaces and for $K\not=0$, we extend the result given in  \cite{dwvw} for CMC translation surfaces.

 A second kind of surfaces of our interest are the homothetical surfaces, where we replace the sum $+$  of a translation surface $z=f(x)+g(y)$ by the multiplication operation $\cdot$ (\cite{wo1}).

 \begin{definition} A homothetical surface $S$ in Euclidean space $\r^3$ is a surface that is a graph of a function $z=f(x)g(y)$, where $f:I\subset\r\rightarrow \r$ and $g:J\subset\r\rightarrow\r$ are two smooth functions.
\end{definition}

 As far as the authors know, the first approach to this kind of surfaces appeared in \cite{wo1} studying the problem of finding minimal homothetical non degenerate surfaces in Lorentz-Minkowski space $\l^3$ (see also \cite{wo2}). Some authors have considered minimal homothetical hypersurfaces in Euclidean   and in semi-Euclidean spaces (\cite{js,wo2}. Our first result concerns minimal surfaces.  Van de Woestyne proved in \cite{wo1} that
the only minimal homothetical non degenerate surfaces in $\l^3$ are planes and helicoids. At the end of \cite{wo1} the author asserted that, up to small changes in the proof, a similar result can be obtained in Euclidean space. In the present paper we do a different  proof of the Euclidean version and in section \ref{s-heli} we prove:

\begin{theorem}\label{t-helicoid}
Planes and helicoids are the only minimal homothetical surfaces in Euclidean space.
 \end{theorem}
The parametrization of helicoid  is not the usual as a ruled surface with base a helix, but as
\begin{equation}\label{pheli}
z(x,y)=(x+b)\tan(c y+d),
\end{equation}
where $b,c,d\in\r$, $c\not=0$.

The third result considers homothetical surfaces in Euclidean space with constant Gauss curvature, obtaining  a complete classification.

\begin{theorem} \label{Ghomo} Let $S$ be a  homothetical surface in Euclidean space with constant Gauss curvature $K$. Then $K=0$. Furthermore, the surface is a cylindrical surface over a plane curve or a surface whose parametrizations is:
\begin{enumerate}
\item $z(x,y)=a e^{bx+cy}$, $a,b,c>0$, or
\item $$z(x,y)=\left(\frac{bx}{m}+d\right)^{m} \left(\frac{cy}{m-1}+e\right)^{1-m},$$
for $b,c,d,e,m\in\r$, $b,c\not=0$, $m\not= 0,1$.
\end{enumerate}
\end{theorem}
This result is proved in sections \ref{sg1} and \ref{sg2} when $K=0$ and $K\not=0$, respectively. Finally, in the last section \ref{s-lore} we extend Theorems \ref{tgt} and \ref{Ghomo} in Lorentz-Minkowski space, obtaining similar results.

\section{Proof of Theorem \ref{tgt}}

Assume $S$ is the sum of the curves $\alpha(s)$ and $\beta(t)$. Locally, $\alpha$ and $\beta$ are graphs on the axis coordinates of $\r^3$ so we may assume that $\alpha(s)=(s,f_1(s),f_2(s))$ and $\beta(t)=(g_1(t),t,g_2(t))$, $s\in I$, $t\in J$, for certain functions $f_1,f_2,g_1,g_2$. Let us observe that if we replace the functions $f_i$ or $g_i$ by an additive constant, the surface changes by a translation of Euclidean space and thus, in what follows, we will take these functions  up to  additive constants. The Gauss curvature in local coordinates $X=X(s,t)$ writes as
$$K=\frac{ln-m^2}{EG-F^2},$$
where $\{E,F,G\}$ and $\{l,m,n\}$ are the coefficients of the first and second fundamental form with respect to $X$, respectively. In our case,  the parametrization of $S$ is $X(s,t)=\alpha(s)+\beta(t)$ and as $\partial^2_{st}X=0$, then $m=0$. The computation of $K$ leads to
\begin{equation}\label{tk}
K=\frac{\left(f_2''-f_1''g_2'+g_1'(f_1''f_2'-f_1'f_2'')\right)\left(g_2''-f_2'g_1''+
f_1'(g_1''g_2'-g_1'g_2'')\right)}{\left((1+f_1'^2+f_2'^2)(1+g_1'^2+g_2'^2)-(f_1'+g_1'+f_2'g_2')^2\right)^2}.
\end{equation}

 \begin{center}{\it Case $K=0$}\end{center}

 Then $l=0$ or $n=0$. Assume $l=0$ and the argument is similar if $n=0$. Thus
\begin{equation}\label{tk1}
f_2''-f_1''g_2'+g_1'(f_1''f_2'-f_1'f_2'')=0
\end{equation}
We distinguish cases.
\begin{enumerate}
\item Assume $f_1''=0$. Then $f_1(s)=as$, $a\in\r$, and \eqref{tk1} gives $f_2''(1-ag_1')=0$.
If $f_2''=0$, then $f_2$ is linear, proving that the curve $\alpha$ is a straight-line and the surface is a cylindrical surface whose base curve is $\beta$. If $f_2''\not=0$, then $a\not=0$. Solving for $g_1$, we obtain $g_1(t)=t/a$. Then $X(s,t)=(s+t/a,as+t,f_2(s)+g_2(t))$ and the surface is the plane of equation $ax-y=0$.
\item Assume $f_1''\not=0$ and $g_1''=0$. Then $g_1(t)=at$, $a\in\r$, and \eqref{tk1} implies
\begin{equation}\label{tgt-ff}
\frac{f_2''+a(f_1''f_2'-f_1'f_2'')}{f_1''}=g_2'.
\end{equation}
As we have a function depending on $s$ equal to a function depending on $t$, then both ones must be a constant $b\in\r$. In particular, $g_2(t)=bt$. Now the curve $\beta$ is a straight-line and the surface is a cylindrical surface with base the curve $\alpha$. Let us notice that under these conditions, equation \eqref{tgt-ff} does not add further information on the curve $\alpha$.
\item Assume $f_1''g_1''\not=0$. Differentiating \eqref{tk1} with respect to $t$, we have
$-f_1''g_2''+g_1''(f_1''f_2'-f_1'f_2'')=0$.  A similar argument as above proves that there exists $a\in\r$ such that
$$\frac{f_1''f_2'-f_1'f_2''}{f_1''}=a=\frac{g_2''}{g_1''}.$$
The identity $g_2''=ag_1''$ implies that $\mbox{det}(\beta',\beta'',\beta''')=0$ and this means that the torsion of $\beta$ is $0$ identically. This proves that $\beta$ is a plane curve. Now we come back to the beginning of the proof assuming that $\beta$ is included in the $yz$-plane (or equivalently, $g_1=0$). We compute $K$ again  obtaining
$$g_2''(f_2''-f_1''g_2')=0.$$
If $g_2''=0$, then $g_2$ is linear and $\beta$ is a straight-line, proving that $S$ is a cylindrical surface with base the curve $\alpha$. If $g_2''\not=0$, then $f_2''-f_1''g_2'=0$ and it follows that there exists $a\in\r$ such that
$$\frac{f_2''}{f_1''}=a=g_2',$$
and so $g_2''=0$, a contradiction.
\end{enumerate}

 \begin{center}{\it Case $K\not=0$}\end{center}

We work following the same ideas as in  \cite{dwvw} by distinguishing two cases: first, we assume that both curves are planar, and second, that only one is planar.
\begin{enumerate}
\item Case that $\alpha$ and $\beta$ are planar. By the result of Liu in \cite{li} , we only consider the case that the curves $\alpha$ and $\beta$ can not lie in planes mutually orthogonal. Let us notice that if the curves lie in parallel planes, the translation surface is (part) of a plane.   Without loss of generality we can assume that $\alpha$ lies in the $xz$-plane and $\beta$ in the plane of equation $x\cos\theta-y\sin\theta=0$, with $\cos\theta,\sin\theta\not=0$. Then
$$\alpha(s)=(s,0,f(s)), \ \beta(t)=(t\sin\theta, t\cos\theta,g(t))$$
with $f$ and $g$ smooth functions on $s$ and $t$, respectively.    The computation of $K$ leads to
$$K=\frac{\cos^2\theta f''g''}{(f'^2+g'^2+\cos^2\theta-2\sin\theta f'g')^2}.$$
Notice that being $K\not=0$, this implies $f'',g''\not=0$. Differentiating with respect to $s$ and with respect to $t$, we obtain respectively
$$\cos^2\theta f'''g''=4K(f'^2+g'^2+\cos^2\theta-2\sin\theta f'g')(f'f''-\sin\theta f''g')$$
$$\cos^2\theta f''g'''=4K(f'^2+g'^2+\cos^2\theta-2\sin\theta f'g')(g'g''-\sin\theta f'g'').$$
Using $f''g''\not=0$, we have
\begin{equation}\label{tkfg}
\frac{f'''}{f''^2}(g'-\sin\theta f')=\frac{g'''}{g''^2}(f'-\sin\theta g').
\end{equation}
Differentiating now with respect to $s$ and next, with respect to $t$, we get
$$\left(\frac{f'''}{f''^2}\right)'g''=f''\left(\frac{g'''}{g''^2}\right)'.$$
Dividing by $f''g''$, we have an identity of two functions, one depending on $s$ and the other one depending on $t$. Then both functions are equal to a constant $a\in\r$. Then there exists $b,c\in\r$ such that
$$\frac{f'''}{f''^2}=af'+b,\ \ \frac{g'''}{g''^2}=ag'+c.$$
Substituting in \eqref{tkfg}, we have
$$a\sin\theta f'^2+b\sin\theta f'+cf'=a\sin\theta g'^2+c\sin\theta g'+bg'.$$
Again we have two functions, one depending on $s$ and other one depending on $t$. Then both functions are constant and hence, $f'$ and $g'$ are constant, in contradiction with   $f''g''\not=0$.
\item Assume that $\alpha$ is a planar curve and $\beta$ does not lie in a plane. After a change of coordinates, we may suppose
$$\alpha(s)=(s,0,f(s)), \ \beta(t)=(g_1(t),t,g_2(t)),$$
for smooth functions $f$, $g_1$ and $g_2$. The contradiction will arrive proving that $\beta$ is a planar curve. For this reason, let us observe that $\beta$ is planar if and only if its torsion vanishes for all $s$, that is, $\mbox{det}(\beta'(t),\beta''(t),\beta'''(t))=0$ for all $t$, or equivalently,
\begin{equation}\label{planarg}
g_1'''g_2''-g_1''g_2'''=0.
\end{equation}
We compute $K$ obtaining
\begin{equation}\label{tkk}
K=\frac{f''(g_2''-f'g_1'')}{(1+g_2'^2+f'^2+f'^2g_1'^2-2f'g_1'g_2')^2}.
\end{equation}
Since $K\not=0$, then $f''\not=0$. In \eqref{tkk},  we move $f''$ to the left hand-side obtaining a function depending only on the variable $s$, the derivative of the right hand-side with respect to $t$ is $0$. This means
$$4(f'g_1'-g_2')(f'g_1''-g_2'')^2+\left(1+f'^2(1+g_1'^2)+2(f'g_1'''-g_2''')(f'g_1'g_2'+g_2'^2)\right)=0.$$
For each fix $t$, we can view this expression as a polynomial equation on $f'(s)$ and thus, all coefficients vanish. Exactly the above equation writes as $\sum_{n=0}^3A_n(t)f'(s)^n=0$. The computations of $A_n$ give:
$$A_0=(1+g_2'^2)g_2'''-4g_2'g_2''^2,\ A_3=-(1+g_1'^2)g_1'''+4g_1'g_1''^2$$
$$A_1=8g_1''g_2'g_2''+4g_1'g_2''^2-(1+g_2'^2)g_1'''-2g_1'g_2'g_2'''$$
$$A_2=-8g_1'g_1''g_2''-4g_1''^2g_2'+2g_1'g_2'g_1'''+(1+g_1'^2)g_2'''.$$
From $A_0=0$ and $A_3=0$ we get
\begin{equation}\label{fi1}
(1+g_i'^2)g_i'''-4g_i'g_i''^2=0, \ \ i=1,2
\end{equation}
obtaining a first integration
\begin{equation}\label{fi2}
g_i''=\lambda_i(1+g_i'^2)^2, \  \lambda_i>0, i=1,2.
\end{equation}
In particular, from \eqref{fi1},
$$g_i'''=4\lambda_i^2 g_i'(1+g_i'^2)^3.$$

Before continuing with the information obtained up here, we write the fact that $\beta$ is a planar curve given in \eqref{planarg} in terms of $g_1'$ and $g_2'$. Then $\beta$ is planar if and only if
\begin{equation}\label{planarg2}
 \lambda_1 g_1'(1+g_1'^2)-\lambda_2 g_2'(1+g_2'^2)=0.
 \end{equation}
From the data obtained for $g_i''$ and $g_i'''$, we now substitute into the coefficients $A_1$ and $A_2$. After some manipulations, the identity $A_1g_2'-A_2g_1'=0$ simplifies into
$$\left[(\lambda_1(1+g_1'^2)g_2'+\lambda_2g_1'(1+g_2'^2)\right]\left[\lambda_1 g_1'(1+g_1'^2)-\lambda_2 g_2'(1+g_2'^2)\right]=0.$$
If the right bracket is zero, then $\beta$ is planar by \eqref{planarg2}, obtaining a contradiction. If the first bracket vanishes, then
$$1+g_2'^2=-\frac{\lambda_1}{\lambda_2}\frac{g_2'}{g_1'}(1+g_1'^2).$$
With this information and \eqref{fi2}, we place it into the coefficient $A_1$, and equation $A_1=0$ is equivalent to the identity
$$ g_1'^4+g_2'^4+g_1'^2+g_2'^2+2g_1'^2g_2'^2=0.$$
Then $g_1'=g_2'=0$, that is, the curve $\beta$ is planar, obtaining a contradiction again. This finishes the proof of Theorem \ref{tgt} for the case $K\not=0$.
\end{enumerate}
\section{Proof of Theorem \ref{t-helicoid}}\label{s-heli}

Assume that $S$ is a homothetical surface which is the graph of $z=f(x)g(y)$ and let $X(x,y)=(x,y,f(x)g(y))$ be a parametrization of $S$. The computation of $H=0$ leads to
\begin{equation}\label{p1}
f''g(1+f^2g'^2)-2ff'^2gg'^2+fg''(1+f'^2g^2)=0.
\end{equation}
Since the roles of $f$ and $g$ in \eqref{p1} are symmetric, we only discuss the cases according to the function $f$. We distinguish cases.
\begin{enumerate}
\item Case $f'=0$. Then  \eqref{p1} gives $fg''=0$. If $f=0$, $S$ is the horizontal plane of equation $z=0$. If $g''=0$, then $g(y)=ay+b$, $a,b\in\r$ and $X(x,y)$ parametrizes the plane of equation $\lambda a y-z=\lambda b$, for some $\lambda\in\r$.
\item Case $f''=0$, $f'\not=0$, and by symmetry, $g'\not=0$. Then $f(x)=ax+b$, for  $a,b\in\r$, $a\not=0$. Now \eqref{p1} reduces into
$$ a^2gg'^2+ g''(1+a^2g^2)=0.$$
Then
$$\frac{g''}{g'}=a^2\frac{gg'}{1+a^2g^2}.$$
 By integrating, there exists a constant $k>0$ such that
$$g'=k(1+a^2g^2).$$
The solution is
$$g(y)=\frac{1}{a}\tan(aky+d), \ d\in\r.$$
 It only rests to conclude that we obtain a helicoid. In such a case, the parametrization of $S$ is
$$X(x,y)=(x,y,f(x)g(y))=(0,y,bg(y))+x(1,0,a g(y))$$
which indicates that the surface is ruled. But it is well known that the only ruled minimal surfaces in $\r^3$ are planes and helicoids and since $g$ is not a constant function, $S$ must be a helicoid (\cite{bc}).
\item Case $f''\not=0$. We will prove that this case is not possible. By symmetry in the discussion of the case, we also suppose $g''\not=0$.  If we divide \eqref{p1} by$ff'^2gg'^2$, we have
$$\frac{f''}{ff'^2g'^2}+\frac{f''f}{f'^2}-2+\frac{g''}{f'^2gg'^2}+\frac{gg''}{g'^2}=0$$
Derive once with respect to $x$ and next with respect to $y$, obtaining
\begin{equation}\label{p2}
 \left( \frac{f''}{ff'^2} \right)' \left( \frac{1}{g'^2} \right)' + \left( \frac{1}{f'^2}\right)' \left( \frac{g''}{gg'^2} \right)' =0.
 \end{equation}
Since $f''g''\not=0$, we divide \eqref{p2} by $(1/g'^2)'(1/f'^2)'$ and we conclude that there exists a constant $a\in\r$ such that
$$ \left( \frac{f''}{ff'^2} \right)' \frac{1}{ \left( \frac{1}{f'^2}\right)' } =a=-\left( \frac{g''}{gg'^2} \right)' \frac{1}{\left( \frac{1}{g'^2} \right)'} .$$
Hence there are constants $b,c\in\r$ such that
$$\frac{f''}{ff'^2} = a\frac{1}{f'^2}+b,\ \ -\frac{g''}{gg'^2} = a\frac{1}{g'^2}+c,$$
or equivalently,
\begin{equation}\label{p3}
f''= f(a+b f'^2),\ \ g''=-g(a+cg'^2).
\end{equation}
Taking into account \eqref{p3}, we replace  $f''$ and $g''$ in \eqref{p1}, obtaining
$$fg\left((a+b f'^2)(1+f^2g'^2)-2f'^2g'^2-(a+cg'^2)(1+f'^2g^2)\right)=0.$$
If we simplify by $fg$ and then we divide by $f'^2g'^2$, we get
$$\frac{c-af^2}{f'^2}+2-b f^2=\frac{b-ag^2}{g'^2}-cg^2.$$
Again we use that each side of  this equation depends on $x$ and $y$, hence there exists $\lambda\in\r$ such that
\begin{equation}\label{p4}
f'^2 = \frac{c-af^2}{\lambda+bf^2-2},\ \ g'^2 =\frac{b-ag^2}{\lambda+cg^2}.
\end{equation}
Differentiating with respect to $x$ and $y$, respectively, we have
\begin{equation}\label{p5}
f''=-\frac{f(bc+a(\lambda-2))}{(\lambda+bf^2-2)^2},\ \ g''=-\frac{ a\lambda+bc}{(\lambda+cg^2)^2}.
\end{equation}
 Let us compare these expressions of $f''$ and $g''$ with those ones that appeared in \eqref{p3} and replace the value of $f'^2$ and $g'^2$ obtained in \eqref{p4}. After some manipulations, we get
$$\left(bc+a(\lambda-2)\right)(\lambda-1+b f^2)=0.$$
$$\left(bc+a \lambda\right)(\lambda-1+cg^2)=0. $$
We discuss all possibilities.
\begin{enumerate}
\item If $bc+a(\lambda-2)=bc+a\lambda=0$, then $a=0$ and $bc=0$. Then \eqref{p5} gives $f''=0$ or $g''=0$, a contradiction.
\item If $bc+a(\lambda-2)=0$ and $c=\lambda-1=0$, we obtain $a=0$. From \eqref{p5}, we get $g''=0$, a contradiction.
\item If $bc+a\lambda=0$ and $b=\lambda-1=0$, then $a=0$ and \eqref{p5} gives $f''=0$, a contradiction.
\item If $b=c=0$ and $\lambda=1$, from the expressions of $f'^2$ and $g'^2$ in \eqref{p4}, we deduce $f'^2=af^2$ and
$g'^2=-ag^2$, that is, $a=0$. Then \eqref{p5} gives $f'=g'=0$, a contradiction again.
\end{enumerate}
\end{enumerate}

\section{Proof of Theorem \ref{Ghomo}: case $K=0$}\label{sg1}

Assume $K=0$. The computation of $K$ gives
\begin{equation}\label{hk}
K=\frac{fgf''g''-f'^2g'^2}{(1+f'^2g^2+f^2g'^2)^2}.
\end{equation}
Hence that if $K=0$, then
\begin{equation}\label{k1}
ff''gg''=f'^2g'^2.
\end{equation}
Since the roles of the function $f$ and $g$ are symmetric in \eqref{k1}, we discuss the different cases according the function $f$.
\begin{enumerate}
\item Case $f'=0$. Then $f$ is a constant function $f(x)=x_0$ and the parametrization of the surface writes as
$X(x,y)=(0,y,x_0g(y))+x(1,0,0)$. This means that $S$ is a cylindrical surface with base curve  a plane curve in the $yz$-plane and the rulings are parallel to the $x$-axis.
\item Case  $f''=0$ and $f',g'\not=0$. Now $f(x)=ax+b$, $a,b\in\r$, $a\not=0$.  Moreover, \eqref{k1} gives $g'=0$ and $g(y)=y_0$ is a constant function.  Now $S$ is the plane of equation $z=x_0(ax+b)$.
\item Case $f''\not=0$. By  the symmetry on the arguments, we also suppose $g''\not=0$.  Equation \eqref{k1} writes as
$$\frac{ff''}{f'^2}=\frac{g'^2}{gg''}$$
As in each side of this equation we have a function depending on $x$ and other depending on $y$, there exists $a\in\r$, $a\not=0$, such that
$$\frac{ff''}{f'^2}=a=\frac{g'^2}{gg''}.$$
A simple integration implies that there exist $b,c>0$ such that
$$f'=bf^a,\ \ g'=cg^{1/a}.$$
\begin{enumerate}
\item Case $a=1$. Then
$$f(x)=pe^{bx}, g(y)=qe^{cy},\ \ p,q>0.$$
\item Case $a\not=1$. Then
$$f(x)=\left((1-a)bx+p\right)^{\frac{1}{1-a}},\ \ g(y)=\left(\frac{a-1}{a}cy+q\right)^{\frac{a}{a-1}},$$
for $p,q\in\r$. This concludes the case $K=0$.
\end{enumerate}
\end{enumerate}

\section{Proof of Theorem \ref{Ghomo}: case $K\not=0$}\label{sg2}

Assume $K$ is a non zero constant. By contradiction, we assume that there exists a homothetical surface $S$ with constant Gauss curvature $K$. Let us observe the symmetry of the expression \eqref{hk} on $f$ and $g$. If $f=0$ or $f'=0$, then \eqref{hk} implies $K=0$, which  is not possible. If $f''=0$, then $f(x)=ax+b$, for some constants $a,b$, $a\not=0$. Then \eqref{hk} writes as
 $$K(1+a^2g^2+(ax+b)^2g'^2)^2+a^2g'^2=0.$$
 This is a polynomial equation on $x$ of degree $4$ because $K\not=0$. Then the leader coefficient, namely, $Ka^4g'^4$, must vanish. This means $g'=0$ and \eqref{hk} gives now $K=0$: contradiction.

Thus  $f''\not=0$. Interchanging the argument with $g$, we also suppose $g''\not=0$. In particular, $fgf'g'\not=0$. We write \eqref{hk} as
$$K=\frac{N}{D^2},\ N=fgf''g''-f'^2g'^2, \ D=1+f'^2g^2+f^2g'^2.$$
By differentiating separately with respect to $x$ and with respect to $y$, and using that $K$ is constant and $N\not=0$,we have
\begin{equation}\label{nd}
\frac{\partial_x N}{\partial_x D}=\frac{\partial_y N}{\partial_y D}.
\end{equation}
Hence we obtain $\partial_x N {\partial_y D-\partial_x D \partial_y N}$. Dividing by $fgf'g'f''^2g''^2$, we get
\begin{eqnarray}\label{p9}
& & -\frac{2f'^2g'^2}{ff''g''^2}+\frac{2f'^2g'^2}{f''^2gg''}+\frac{3f'^2g}{ff''g''}-\frac{3fg'^2}{f''gg''}+
\frac{f'f'''g}{f''^2g''}-\frac{fg'g'''}{f''g''^2}\nonumber \\
& &+ \mbox{ function on $x$}+\mbox{function on $y$}=0.
\end{eqnarray}
Derive \eqref{p9} with respect to $x$ and next with respect to $y$, obtaining
\begin{equation}\label{p92}
-\left(\frac{2f'^2} {ff''}\right)'\left(\frac{g'^2}{g''^2}\right)'+\left(\frac{2f'^2}{f''^2}\right)'\left(\frac{g'^2}{gg''}\right)'+
\left(\frac{3f'^2}{ff''}\right)'\left(\frac{g}{g''}\right)'-\left(\frac{3f}{f''}\right)'\left(\frac{g'^2}{gg''}\right)'=0.
\end{equation}
In this expression, we want to divide by
\begin{equation}\label{p93}
\left(\frac{f}{f''}\right)'\left(\frac{g}{g''}\right)'.
\end{equation}
In order to do it, we must assure that this function is not zero. We stop here the proof of Theorem \ref{Ghomo}, and let us discuss this possibility.

{\it Claim.} The case $\left(\frac{f}{f''}\right)'\left(\frac{g}{g''}\right)'=0$ is not possible.

\begin{proof} (of the claim)
Without loss of generality, we suppose that $\left(\frac{f}{f''}\right)'=0$. Then there exists $\lambda\not=0$ such that $f''=\lambda f$. The solution of this equation depends on the sign of $\lambda$.  We may reduce the discussion to the cases $\lambda=1$ and $\lambda=-1$. This is possible by a change of coordinates $(x,y)\longmapsto k(x,y)$ on the initial parametrization $X(x,y)$ of $S$, obtaining a new surface which is homothetical to the original and this does not change the fact that $K$ is a non zero constant.
\begin{enumerate}
\item Case $\lambda=1$. Then the solution of $f''=f$ is
$f(x)=a\cosh(x)+b\sinh(x)$ with $a,b\in\r$. We substitute in \eqref{hk}, and after some manipulations, we obtain an expression
\begin{equation}\label{hkcs}
\sum_{n=0}^4A_n(y)\cosh(nx)+B_n(y)\sinh(nx)=0.
\end{equation}
 Since the functions $\{\cosh(nx),\sinh(nx):n\in{\mathbb N}\}$ are linearly independent, then $A_n=B_n=0$ for $0\leq n\leq 4$. For $n=4$, we have that $A_4=B_4=0$ writes as
$$(a^4+6a^2b^2+b^4)(g^2+g'^2)=0,\ (a^2+b^2)(g'^2+g^2)=0,$$
obtaining a contradiction.
\item Case  $\lambda=-1$. The solution  of $f''=-f$ is $f(x)=a\cos(x)+b\sin(x)$ with $a,b\in\r$. Now \eqref{hk} is
$$\sum_{n=0}^4A_n(y)\cos(nx)+B_n(y)\sin(nx)=0.$$
Again all coefficients $A_n$ and $B_n$ must vanish. In particular, from $A_4=B_4=0$, we obtain
$$K(a^4-6a^2b^2+b^4)(g'^2-g^2)=0,\ abK(a^2-b^2)(g'^2-g^2)=0.$$
If $g'^2\not=g^2$, then we deduce from both equations $a=b=0$, a contradiction. Thus $g'^2=g^2$, that is, $g(y)=ce^y$ or $g=ce^{-y}$, $c>0$. We do the argument for  $g(y)=ce^y$ (the case $g(y)=ce^{-y}$ is analogous). Then \eqref{hk} reduces into
$$(a^2+b^2)^2c^4 K e^{4y}+(a^2+b^2)c^2(1+2K)e^{2y}+K=0,$$
which is a polynomial equation on $e^{2y}$. Then the coefficients corresponding to $e^{2y}$ and $e^{4y}$ vanish, obtaining a contradiction.
\end{enumerate}
\end{proof}

Once proved the claim, we follow with the proof of Theorem \ref{Ghomo}. First, we introduce the notation
$$A=\frac{f}{f''},\  B=\frac{g}{g''}.$$
Divide \eqref{p92} by $A'B'$. By differentiating with respect to $x$ and next with respect to $y$, we conclude that there exists $a\in\r$ such that
$$\left(\frac{\left(\frac{f'^2} {ff''}\right)'}{A'}\right)'=a\left(\frac{\left(\frac{f'^2}{f''^2}\right)'}{A'}\right)'$$
$$\left(\frac{\left(\frac{g'^2} {gg''}\right)'}{B'}\right)'=a\left(\frac{\left(\frac{g'^2}{g''^2}\right)'}{B'}\right)'.$$
Then there are $b_1,b_2,c_1,c_2\in\r$ such that
\begin{eqnarray}
\frac{f'^2} {ff''}&=&a\frac{f'^2}{f''^2}+b_1 A+b_2\label{p9f}\\
\frac{g'^2} {gg''}&=&a\frac{g'^2}{g''^2}+c_1 B+c_2\label{p9g}
\end{eqnarray}
Substituting \eqref{p9f} and \eqref{p9g} into \eqref{p92}, we get
\begin{equation}\label{p94}
(3a-2c_1)B'\left(\frac{f'^2}{f''^2}\right)'+3(b_1-c_1)A'B'+
(3a+2b_1)A'\left(\frac{g'^2}{g''^2}\right)'=0.
\end{equation}
Our argument depends on the coefficients of the above equation and we will take into account the symmetry on the reasoning changing $f$ by $g$.  ,
\begin{center}Case $a=0$\end{center}
\begin{enumerate}
\item Sub-case $b_1=0$. From \eqref{p9f},  $ff''=\lambda f'^2$ for $\lambda\not=0$.   The solutions are
$f(x)=n e^{mx}$, $n,m>0$ ($\lambda=1$) or $f(x)=(1-\lambda) (mx+n)^{\frac{1}{1-\lambda}}$ ($\lambda\not=1$), $m>0$.
In the first case, \eqref{hk} writes as a polynomial equation on $e^{mx}$. Then the coefficients must vanish, but the leader coefficient corresponding to  $e^{4mx}$ is $Kn^4(m^2g^2+g'^2)$, a contradiction. In the second case, \eqref{hk} is a polynomial equation on $x$ and the coefficient for $x^{\frac{4}{1-\lambda}}$ is $K(1-\lambda)^4m^{\frac{4}{1-\lambda}}g'^4$. Because this must vanish, then $g'=0$, a contradiction.
\item Sub-case $b_1\not=0$, and by symmetry, $c_1\not=0$. Then \eqref{p94} writes as
$$-2c_1B'\left(\frac{f'^2}{f''^2}\right)'+3(b_1-c_1)A'B'+2b_1 A'\left(\frac{g'^2}{g''^2}\right)'=0.$$
Divide by $A'B'$. Then
$$\frac{-2c_1 \left(\frac{f'^2}{f''^2}\right)'}{A'}+3(b_1-c_1)=-2b_1\frac{\left(\frac{g'^2}{g''^2}\right)'}{B'}.$$
Again we have two functions, one depending on $x$ and the other one depending on $y$. This means that there exists $\lambda,m,n\in\r$ such that
$$ \frac{f'^2}{f''^2}=\frac{\lambda+3(b_1-c_1)}{2c_1}A+m,\ \ \frac{g'^2}{g''^2}=\frac{\lambda}{2b_1}B+n.$$

Dividing both expressions by the ones given in \eqref{p9f} and \eqref{p9g} and taking into account that $a=0$, we obtain
$$A=\frac{\displaystyle\frac{\lambda+3(b_1-c_1)}{2c_1}A+m}{b_1 A+b_2},\ \ B=\frac{\displaystyle\frac{\lambda}{2b_1}B+n}{c_1 B+c_2},$$
or equivalently
$$b_1A^2+\left(b_2-\frac{\lambda+3(b_1-c_1)}{2c_1}\right)A-m=0$$
$$c_1 B^2+\left(c_2-\frac{\lambda}{2b_1}\right)B-n=0.$$
By solving both expressions for $A$ and $B$, we obtain that $A$ and $B$ are constant, a contradiction with the Claim.  Thus the only possibility is that the coefficients of the above polynomials on $A$ and on $B$ must vanish. This implies $b_1=c_1=0$, a contradiction.
\end{enumerate}
\begin{center}Case $a\not=0$\end{center}
Let us observe that the three coefficients  in \eqref{p94} can not be all $0$ because this implies $a=0$. In a first discussion, we suppose that one of these coefficients is $0$. By using that $A'B'\not=0$ and $f'g'\not=0$, then necessarily the other two coefficients are not $0$.
 \begin{enumerate}
 \item Sub-case $b_1-c_1=0$. Then
 $$(3a-2b_1)B' \left(\frac{f'^2}{f''^2}\right)'+(3a+2b_1)A'\left(\frac{g'^2}{g''^2}\right)'=0.$$
 Dividing by $A'B'$, there exists $\lambda,m,n\in\r$ such that
\begin{equation}\label{ay1}
\frac{f'^2}{f''^2}=\frac{\lambda}{3a-2b_1}A+m,\ \ \frac{g'^2}{g''^2}=\frac{-\lambda}{3a+2b_1}B+n.
\end{equation}
 We substitute in    \eqref{p9f} and \eqref{p9g}, obtaining
 $$\frac{f'^2} {ff''}=a\left(\frac{\lambda}{3a-2b_1}A+m\right) +b_1 A+b_2$$
$$\frac{g'^2} {gg''}=a\left(\frac{-\lambda}{3a+2b_1}B+n\right)+b_1 B+c_2$$
Dividing again these expressions by \eqref{ay1}, we have
$$(a\alpha+b_1)A^2+(am+b_2-\alpha)A-m=0,$$
$$(a\beta+b_1)B^2+(an+c_2-\beta)B-n=0,$$
where
$$\alpha=\frac{\lambda}{3a-2b_1}, \ \beta=\frac{-\lambda}{3a+2b_1}.$$
Again, the above expressions prove that, either $A$ and $B$ are constant or the coefficients of both polynomials on $A$ and on $B$ vanish. Since the first case is impossible, then  $n=m=0$ and
$$a\alpha+b_1=0, b_2-\alpha=0, a\beta+b_1=0,c_2-\beta=0.$$
Then $\alpha=\beta$ and this implies $a=0$: contradiction.
 \item Sub-case $3a-2c_1=0$. From \eqref{p94}, there exists $m\in\r$ such that
 \begin{equation}\label{ay2}
 \frac{g'^2}{g''^2}=-3\frac{b_1-c_1}{3a+2b_1}B+m
 \end{equation}
 and substituting in \eqref{p9g},
 $$\frac{g'^2}{gg''}=a\left(-3\frac{b_1-c_1}{3a+2b_1}B+m\right)+c_1B+c_2.$$
 Dividing by \eqref{ay2}, we have
 $$B=\frac{\alpha B+m}{(a\alpha+c_1)B+am+c_2},\ \alpha=-3\frac{b_1-c_1}{3a+2b_1}.$$
 Then
 $$(a\alpha+c_1)B^2+(am+c_2-\alpha)B-m=0.$$
 Since $B$ is not constant, then $m=0$ and $\alpha=c_2=-c_1/a$. By using $3a-2c_1=0$, we get $a=0$: contradiction.
 \item Sub-case $3a+2b_1=0$. It is analogous than the sub-case $3a-2a_1=0$ by changing the roles of $f$ and $g$.
 \item Sub-case $(3a-2c_1)(b_1-c_1)(3a+2b_1)\not=0$. Dividing by $A'B'$ in  \eqref{p94} we have
$$ (3a-2c_1)\frac{\left(\frac{f'^2}{f''^2}\right)'}{A'}+3(b_1-c_1)=-(3a+2b_1)\frac{\left(\frac{g'^2}{g''^2}\right)'}{B'}=\lambda$$
for some $\lambda\not=0$. There exist $m,n\in\r$ such that
\begin{equation}\label{ay3}
\frac{f'^2}{f''^2}=\frac{\lambda-3(b_1-c_1)}{3a-2c_1}A+m,\ \frac{g'^2}{g''^2}=\frac{-\lambda}{3a+2b_1}B+n.
\end{equation}
We substitute in    \eqref{p9f} and \eqref{p9g}, obtaining
$$\frac{f'^2} {ff''}=a\left(\alpha A+m\right) +b_1 A+b_2$$
$$\frac{g'^2} {gg''}=a\left(\beta B+n\right)+c_1 B+c_2$$
where
$$\alpha=\frac{\lambda-3(b_1-c_1)}{3a-2c_1}, \ \beta=\frac{-\lambda}{3a+2b_1}.$$
Dividing again these expression by \eqref{ay3}, we have
$$(a\alpha+b_1)A^2+(am+b_2-\alpha)A-m=0,$$
$$(a\beta+c_1)B^2+(an+c_2-\beta)B-n=0.$$
As $A$ and $B$ are not constant,  we get $m=n=0$ and $a\alpha+b_1=b_2-\alpha=a\beta+c_1=c_2-\beta=0$. Again we conclude $a=0$: contradiction.
 \end{enumerate}

\section{The Lorentzian case}\label{s-lore}

We consider the Lorentzian-Minkowski space $\l^3$, that is, $\r^3$ endowed with the metric $(dx)^2+(dy)^2-(dz)^2$. A surface immersed in $\l^3$ is said non degenerate if the induced metric on $S$ is non degenerate. The induced metric can only be of two types: positive definite and the surface is called spacelike, or a Lorentzian metric, and the surface is called timelike. For both types of surfaces, it is defined the mean curvature $H$ and the Gauss curvature, which have the following expressions in local coordinates $X=X(s,t)$:
$$H=\epsilon\ \frac12\frac{lG-2mF+nE}{EG-F^2},\ \ K= \epsilon\ \frac{ln-m^2}{EG-F^2},$$
where $\epsilon=-1$ is $S$ is spacelike and $\epsilon=1$ if $S$ is timelike. Here $\{E,F,G\}$ and $\{l,m,n\}$ are the coefficients of the first and second fundamental forms with respect to $X$, respectively. See \cite{lo} for more details. Again we ask for those translation and homothetical surfaces in $\l^3$ with constant mean curvature and constant Gauss curvature. Recall that the property of a surface to be a translation surface or a homothetical surface is not metric but it is given by the affine structure of $\r^3$ and the multiplication of real functions of $\r$.

We generalize the results obtained in the previous sections for non degenerate surfaces of $\l^3$. The proofs are similar, and we omit the details.

\begin{enumerate}

\item Extension of Theorem \ref{tgt}. Assume that $S$ is a translation surface. The computation of $K$ gives
$$K=-\frac{\left(f_2''-f_1''g_2'+g_1'(f_1''f_2'-f_1'f_2'')\right)\left(g_2''-f_2'g_1''+
f_1'(g_1''g_2'-g_1'g_2'')\right)}{\left((1+f_1'^2-f_2'^2)(1+g_1'^2-g_2'^2)-(f_1'+g_1'-f_2'g_2')^2\right)^2}.$$
If $K=0$, then the numerator coincides with the one in \eqref{tk} and the conclusion is that $S$ is a cylindrical surface.
 In the case $K\not=0$, the results asserts that, under the same hyposthesis, there are no further examples. We discuss the cases that $\alpha$ and $\beta$ lies in two non-orthogonal planes and when one curve is planar. In the former case, the expression of $K$ is
 $$K=-\frac{\cos^2\theta f''g''}{(-f'^2-g'^2+\cos^2\theta+2\sin\theta f'g')^2}.$$
 The proof works the same.  In the second case,
$$K=-\frac{f''(g_2''-f'g_1'')}{(1-g_2'^2-f'^2-f'^2g_1'^222f'g_1'g_2')^2}.$$
Again, the proof is identical because we can move $f''$ to the left hand-side, differentiating with respect to $t$ and observing that it appears an expression which is a polynomial on the function $f'$.
\item Extension of Theorem \ref{t-helicoid}. As we have pointed out, the result in Lorenztian case was proved in \cite{wo1}.
\item Extension of Theorem \ref{Ghomo}.  Assume now that $S$ is a homothetical surface and we study those surfaces with constant Gauss curvature. If $S$ is spacelike, then the surface is locally a graph on the $xy$-plane and  $S$ writes as $z=f(x)g(y)$. The expression of $K$ is
$$K=-\frac{fgf''g''-f'^2g'^2}{(1-f'^2g^2-f^2g'^2)^2},\ \mbox{with $1-f'^2g^2-f^2g'^2>0$}.$$
If $S$ is timelike, then the surface is locally a graph on the $xz$-plane or on the $yz$-plane. Without loss of generality, we assume that the surface writes as $x=f(y)g(z)$. Now the Gauss curvature $K$ is
$$K=-\frac{fgf''g''-f'^2g'^2}{(1+f^2g'^2-f'^2g^2)^2},\ \mbox{with $1+f^2g'^2-f'^2g^2<0$}.$$
Because both expression are similar as in \eqref{hk} and the arguments are the same as in Euclidean space, we only give the statements. If $K\not=0$, then there are not exist homothetical (spacelike or timelike) surfaces with constant Gauss curvature $K$. If $K=0$, then $fgf''g''-f'^2g'^2=0$, which it is the same that \eqref{k1}. Then the conclusion is:
 \begin{enumerate}
 \item The surface is  a cylindrical surface whose base curve is a planar curve contained in one of the three coordinates planes and the ruling are orthogonal to this plane.
 \item  The function $z=f(x)g(y)$ agrees as in Theorem \ref{Ghomo}, 1) and 2).
\end{enumerate}
\end{enumerate}



\begin{thebibliography}{11}
\bibitem{bc} J. L. M. Barbosa, A. Colares, A. G.: \emph{Minimal Surfaces in $R^3$},  Lecture Notes in Mathematics, 1195, Springer-Verlag, Berlin, 1986.

\bibitem{da} J. G. Darboux,  \emph{Th\'eorie G\'enerale des Surfaces}, Livre I, Gauthier-Villars, Paris, 1914.

\bibitem{dwvw} F. Dillen, I. Van de Woestyne, L. Verstraelen, J. T. Walrave,
The surface of Scherk in $E^3$: a special case in the class of minimal surfaces defined as the sum of two curves,
\emph{Bull. Inst. Math. Acad. Sin.} 26 (1998), 257--267.

\bibitem{js} L. Jiu, H.  Sun, On minimal homothetical hypersurfaces,
\emph{Colloq. Math.} 109 (2007), 239--249.

\bibitem{li} H. Liu,  Translation surfaces with constant mean curvature in 3-dimensional spaces, \emph{J. Geom.} 64 (1999), 141--149.

\bibitem{lo} R. L\'opez,  Differential Geometry of curves and surfaces in Lorentz-Minkowski space, \emph{Internat.  Electronic J. Geom.}  7 (2014), 44--107.

\bibitem{ni} J. C. C. Nitsche,  \emph{Lectures on Minimal Surfaces}, Cambridge University Press, Cambridge, 1989.

\bibitem{sc}  H. F. Scherk,  Bemerkungen \"{u}ber die kleinste Fl\"{a}che innerhalb gegebener Gren-
zen, \emph{J. Reine Angew. Math.} 13 (1835), 185--208.

\bibitem{wo1} I. Van de Woestyne,  A new characterization of the helicoids, \emph{Geometry and topology of submanifolds, V} (Leuven/Brussels, 1992), 267--273, World Sci. Publ., River Edge, NJ, 1993.

\bibitem{wo2} I. Van de Woestyne,  Minimal homothetical hypersurfaces of a semi-Euclidean space, \emph{Results Math.} 27 (1995), 333--342.

 \end{thebibliography}
\end{document}